\documentclass{amsart}
\usepackage{amsfonts,amscd}
\usepackage[dvips]{graphicx}
\input xy
\xyoption{all}

\theoremstyle{definition}

\numberwithin{equation}{section}

\def\Z{{\mathbb Z}}

\def\R{{\mathbb R}}

\def\P{{\mathbb P}}
\def\gp{\gamma_{p}}

\begin{document}

\title[The Isotopy Problem for Symplectic \textbf{4}-manifolds]{The Isotopy Problem for Symplectic \textbf{4}-manifolds}

\author{Stefano Vidussi}
\address{Department of Mathematics, Kansas State University,
Manhattan, KS 66506, USA}
\email{vidussi@math.ksu.edu} \thanks{The author was supported in part
by NSF grant \#0306074.}


\begin{abstract}
The purpose of this paper is to present some results on the existence of  homologous, nonisotopic symplectic or lagrangian surfaces embedded in a simply connected symplectic $4$-dimensional manifold. 
\end{abstract}

\maketitle

\section{Introduction}

\noindent Let $X$ be a closed simply connected symplectic $4$-manifold and let $\omega \in \Omega^{2}(X)$ be a symplectic form on $X$, i.e. a closed, nondegenerate $2$-form. A closed (oriented) embedded surface $\Sigma \hookrightarrow X$ is called \\ - \textit{symplectic} if $\omega$ restricts to a symplectic form on $\Sigma$, i.e. $\omega|_{\Sigma} > 0$;
\\ - \textit{lagrangian} if $\omega$ restricts to the trivial form on $\Sigma$, i.e. $\omega|_{\Sigma} = 0$.

Given two submanifolds representing the same homology class, we say that they are isotopic if they can be joined by a family of embedded submanifolds. By the Isotopy Extension Theorem (see e.g. \cite{GS}), this is equivalent to the existence of an ambient isotopy of $X$ connecting them. Isotopic submanifolds are essentially the same, from the point of view of the smooth topology of $X$, so in order to study the topological aspects of symplectic and lagrangian submanifolds of $X$ we need to take into account the equivalence relation above. This leads to the following 
\\

\noindent \textbf{Problem} (P): Fix a homology class $\alpha \in H_{2}(X,\Z)$; classify all the symplectic (lagrangian) representatives of $\alpha$ up to the equivalence relation above.  
\\

This problem was first posed, in the study of lagrangian $\R^{2}$'s linear outside a ball in $\R^{4}$ (where the homology condition becomes irrelevant) by Arnol'd in \cite{A}, and further generalized, also in the lagrangian case, by Eliashberg and Polterovich (see \cite{ElP3}). The interest in the symplectic version of this problem was triggered more recently by Siebert and Tian (see \cite{Ti}).

There are several ways in which we can interpret (P). The first starts from the observation that two submanifolds of a simply connected $4$-manifold are homologous if and only if they are images of homotopic embeddings (this follows from a quite straightforward exercise in algebraic topology). Therefore the study of the representatives of a given homology class, up to isotopy, is a suitable generalization of classical knot theory. Problem (P) amounts therefore to symplectic (lagrangian) \textit{surface-knot theory} for $(M,\alpha)$. The second interpretation stems from the classical observation that if we consider the case when $X$ is a K\"ahler surface, and we look at complex representatives of the same homology class, these are always isotopic (even in a stronger sense, i.e. through complex representatives). Problem (P) can be interpreted therefore as the investigation of the natural counterpart of that property in the symplectic setup, with the added condition that, in a symplectic manifold, we have two kinds of natural submanifolds. (Observe that, while symplectic submanifolds generalize the notion of complex submanifolds, they are defined with an open condition, while lagrangian submanifolds are defined, like the complex ones, with a closed one; in this sense, both kinds of submanifolds share some of the properties of complex ones.) With this interpretation, Problem (P) enters in the larger question of understanding whether the class of symplectic manifolds, that lies somewhere between the class of K\"ahler surfaces and general $4$-manifolds, exhibits features typical of the former or of the latter. To clarify this comment, we want to point out that for a smooth $4$-manifold $X$ there is a rather general construction (see \cite{FS2}) of infinitely many nonisotopic representatives of a homology class $\alpha$, so that the smooth counterpart of problem (P) (i.e. surface-knot theory for $(M,\alpha)$) can be extremely rich, as opposed to ``complex" knot theory for K\"ahler manifolds which is, as noted, trivial.

Problem (P), in the way it is stated, is clearly very ambitious, and what we can realistically hope at the moment is to prove that for a pair $(M,\alpha)$ there is a unique symplectic (lagrangian) representative, or else build infinitely many nonisotopic representatives.

The first results around (P) have all been in the direction of uniqueness. Just to mention few, Eliashberg and Polterovich proved in \cite{ElP2} uniqueness of lagrangian representatives for $\R^{2}$ in $\R^{4}$, answering Arnol'd's question; also, they proved in \cite{ElP1} an analogous result for the class of the section in $T^{*}S^{2}$ and $T^{*}T^{2}$. In the symplectic case, Siebert and Tian proved in \cite{ST} that a symplectic surface in $\P^{2}$ homologous to an algebraic curve of degree $d$ is in fact isotopic to it (at least for $d \leq 17$).

Such results, and the nature of their proofs, might have suggested the possibility of obtaining a uniqueness result holding for all pairs $(M,\alpha)$, as in the K\"ahler case. Instead, Fintushel and Stern have shown in \cite{FS3} that this is not the case, exhibiting the first example of homologous, nonisotopic symplectic tori, in a large class of symplectic manifolds (including, for example, elliptic surfaces). The first examples of nonuniqueness for the lagrangian case have appeared in \cite{V2}.

The aim of this paper is to present and discuss several examples where we can partially answer (P) showing nonuniqueness of 
symplectic or lagrangian representatives for some pair $(M,\alpha)$. All these constructions, in a more or less direct way, start from the ideas contained in \cite{FS3}, and we will try to point out the several common aspects. In Section \ref{const} we are going to discuss the actual construction of infinite families of symplectic (lagrangian) representatives for suitable pairs $(M,\alpha)$, while in Section \ref{noniso} we will discuss how (and to which level) we can distinguish the elements of these families. In these two sections we will restrict the discussion to genus one surfaces. In Section \ref{last}, instead, we are going to show a more recent construction of examples for genus greater than one. 

The results presented below, when not explicitly attributed, can be considered to be a sort of ``common folklore", or can be easily obtained by looking at related papers (of the author or others), quoted in the References. It should be clear, however, that much credit is due to the ideas contained in \cite{FS1} and \cite{FS3}. Finally, for obvious reasons of brevity, this paper does not aspire to be a complete survey on the topic, and follows instead tastes and preferences of the author, even at the price of omitting significant results. 

\section{Construction of the examples}\label{const} 

It is quite clear that describing the symplectic (lagrangian) submanifolds of a symplectic $4$-manifold starting from the definition is a rather difficult task, as the definition itself does not shed too much light on global properties of such submanifolds. The basic idea underlying this section will be to use a class of symplectic $4$-manifolds where we are allowed to reduce the problem, at least in part, to a $3$-dimensional problem, and we can translate the property of being symplectic or lagrangian into a quite simpler problem of curves in a fibered $3$-manifold.

The class of $4$-manifolds that we are going to consider, and that has saved the day in many problems of symplectic topology, is the class of link surgery manifolds. Precisely, we are going to exploit the following dictionary (inspired from \cite{McMT}):
\vspace*{7pt}
\begin{center}
\begin{tabular}{|c|c|}
      \hline   \textbf{3-manifold} & \textbf{4-manifold} \\
      \hline link $L \subset S^{3}$ & link surgery manifold \\
      \hline Alexander polynomial & Seiberg-Witten polynomial \\ 
      \hline fibered link & symplectic $4$-manifold \\
	\hline fibration $1$-form & symplectic $2$-form \\
	\hline curve transverse to fibers & symplectic torus \\
	\hline curve on a fiber & lagrangian torus \\
      \hline
    \end{tabular} \end{center} 
\vspace*{15pt}

To explain the dictionary, we start with the definition of link surgery manifolds (see \cite{FS1}) for a fibered link $L$. 
Consider an $n$-component oriented link $L \subset S^{3}$ whose exterior fibers over $S^{1}$ with fiber $\Sigma$, i.e. $$  S^{3} \setminus \nu L \stackrel{\Sigma}{\longrightarrow} S^{1}. $$ Note that we do not require $\Sigma$ to span the link $L$, namely the (Poincar\'e dual of the) class of $\Sigma$ in $H^{1}(S^{3} \setminus \nu L,\Z) = \Z^{n}$ may be different from $(1,...,1)$. In particular, we will consider the case where $\Sigma$ spans one link component, e.g. $[\Sigma] = (1,0,...,0)$. If $[\Sigma] = (m_{1},...,m_{n})$, the fiber $\Sigma$ intersects the $i$-th boundary component $\partial \nu K_{i}$, up to isotopy, in the curve 
\begin{equation} \label{slope} d_{i} \sigma_{i} = - (\sum_{j \neq i} m_{j} \mbox{lk}(K_{i},K_{j})) 
\mu(K_{i}) + m_{i} \lambda(K_{i}) \end{equation} 
where $\sigma_{i}$ is a simple closed curve and $d_{i} > 0$ is some integer coefficient determined by the equation above (see, e.g., \cite{EN}). Complete $\sigma_{i}$ to a basis $(\rho_{i}, \sigma_{i})$ for $\partial \nu K_{i}$ . For example, if $n=1$, $\sigma = \lambda(K)$ and we can choose $\rho = \mu(K)$.
Next, choose a family $X_{i}$ of symplectic $4$-manifolds, each containing a framed, symplectic, self-intersection $0$ torus $F_{i}$. We can now define the $4$-manifold
\begin{equation} \label{glue} X_{L} = \coprod_{i=1}^{n} (X_{i} \setminus \nu F_{i} ) \cup_{\partial \nu F_{i} = S^{1} \times \partial \nu K_{i}} S^{1} \times
(S^{3} \setminus \nu L). \end{equation} The gluing map on the boundary $3$-tori is defined in such a way to
identify $F_{i}$ with $S^{1} \times \rho_{i}$ and, reversing the orientation, the meridian to the torus $F_{i}$ with $\sigma_{i}$. The manifold above depends on various choices, besides the link $L$ and the pairs $(X_{i},F_{i})$, but we are going to omit such dependence from the notation as this will not affect our discussion. The homotopy type of $X_{L}$ can be easily determined from the $X_{i}$'s and the gluing data. In particular, in the cases we will consider $X_{L}$ will be simply connected. An important feature of this manifold is that
we can explicitly compute the Seiberg-Witten polynomial of $X_{L}$ (and, more generally, of any link surgery manifold) in terms of the relative Seiberg-Witten invariants of $(X_{i},F_{i})$ and the symmetrized multivariable Alexander polynomial of the link $L$. The work of \cite{FS1} and \cite{T} shows that, for $n > 1$ (the case we will be interested in), 
\begin{equation} \label{alex} SW_{X_{L}} = (\coprod_{i=1}^{n} SW_{(X_{i},F_{i})}) \cdot \Delta_{L}(t_{1}^{2},...,t^{2}_{n}). \end{equation}
In this equation $t_{i}$ represents the (Poincar\'e dual to the) homology class of $S^{1} \times \mu(K_{i})$, and we are identifying in the product the homology classes according to the gluing of (\ref{glue}). Remember that 
the relative Seiberg-Witten polynomial of $(X_{i},F_{i})$ is defined as the Seiberg-Witten polynomial of $X_{i} \#_{F_{i} = F} E(1)$, where $E(n)$ is the simply connected elliptic surface without singular fibers of Euler characteristic $12n$ (in particular $SW_{(E(1),F)} = SW_{E(2)} = 1$, and we will be mainly interested in this case).

To show that $X_{L}$ admits a symplectic structure, we need a second presentation of it as symplectic fiber sum of
symplectic $4$-manifolds. In order to do so, perform Dehn surgery on 
$S^{3}$ along $L$, using as slopes on each boundary component the curves $\sigma_{i}$ described in (\ref{slope}): this gives a closed three manifold $N_{L}$. By extending the fibration of $S^{3} \setminus \nu L$ with the fibration of $S^{1} \times D^{2}$ of degree $d_{i}$ on each Dehn filling, $N_{L}$ fibers over $S^{1}$ with fiber $\hat \Sigma$, the closed surface obtained by capping off each boundary component of $\Sigma$ with a disk (contained in the suitable Dehn filling). The cores $C_{i}$ of the Dehn fillings
(i.e. the so called dual link) are transverse, up to isotopy, to ${\hat \Sigma}$. (Note that, unless all divisibilities $d_{i}$ are equal to one, these are not sections of the fibration.) Now we can associate to a fibration of $N_{L}$ a closed nondegenerate $1$-form $\alpha \in \Omega^{1}(N_{L})$, defined up to isotopy. If we denote by $\beta \in \Omega^{2}(N_{L})$ a closed two-form on $N_{L}$ that restricts to a volume form on the fiber ${\hat \Sigma}$, the $4$-manifold 
$S^{1} \times N_{L}$ admits a symplectic $2$-form $\omega = dt \wedge \alpha + \epsilon \beta$, where $\epsilon$ is a sufficiently small constant. Consider the cores $C_{i}$; because of the transversality condition we have $\alpha(C_{i}) > 0$, which implies that \begin{equation} \label{tra} \omega |_{S^{1} \times C_{i}} = dt \wedge \alpha|_{S^{1} \times C_{i}} = \alpha|_{C_{i}} > 0 \end{equation} i.e. the tori $S^{1} \times C_{i}$ 
are symplectic. Moreover, they are naturally endowed with a framing inherited from the framing of the $C_{i}$'s, that we can assume to be the one identified by $\rho_{i}$.
By suitably scaling the symplectic form on the $X_{i}$'s, we can arrange that the symplectic volume on $S^{1} \times C_{i}$ coincides with the one of $F_{i} \subset X_{i}$. We can then write $X_{L}$ as symplectic fiber sum: 
$$ X_{L} = \coprod_{i=1}^{n} X_{i} \#_{F_{i}=S^{1} \times C_{i}}
S^{1} \times N. $$ Gompf's theory of symplectic fiber sums (\cite{G}) guarantees then that $X_{L}$ admits a symplectic form which restricts, outside the gluing locus, to the ones of the factors. 

There is a second result (see \cite{McMT}) that allows us to control some of the homology classes of $X_{L}$, namely the fact that the natural map $$ H_{1}(S^{3} \setminus \nu L,\Z) \rightarrow H_{2}(S^{1} \times (S^{3} \setminus \nu L),\Z) \rightarrow 
H_{2}(X_{L},\Z), $$ obtained by composing K\"unneth isomorphism and inclusion is injective. As a consequence of the previous formula if $\gp$ is a curve in $S^{3} \setminus \nu L$ the homology class of the torus $S^{1} \times \gp$ in $X_{L}$ is given by \begin{equation} \label{hom}
[S^{1} \times \gp ] = \sum_{i=1}^{n} \mbox{lk}(\gp,K_{j}) [S^{1} \times \mu(K_{j})]. \end{equation}

If we apply the construction above choosing as link $L$ the trivial knot $K$ and a symplectic manifold $X$, we obtain of course $X_{K} = X$; the advantage of the presentation (\ref{glue}) is that we have an explicit control of the symplectic topology of the neighborhood of $F$. When $L$ is a nontrivial fibered knot, there is only one choice for $\Sigma$, and $X_{K}$ differs from $X$ by the fact that we have substituted a neighborhood of $F$ with a homology $F \times D^{2}$ that is actually a $\Sigma$-bundle over $T^{2}$, and has a richer topology than the product bundle.

The strategy to obtain families of symplectic tori in $X_{L}$ follows directly from the observations above. Namely, if we can identify a family $\{ \gp \}$ of closed curves in $S^{3} \setminus \nu L$ that are transverse to the fibration, these will define symplectic tori $S^{1} \times \gp$ in $X_{L}$, as the steps leading to (\ref{tra}) follow \textit{verbatim}. Moreover, if the 
$\gp$'s represent the same homology class of $S^{3} \setminus \nu L$, the symplectic tori $S^{1} \times \gp$ will represent the same homology class of $X_{L}$, according to (\ref{hom}). If the curves $\gp$ are not isotopic (in $S^{3} \setminus \nu L$), it is reasonable to expect that the tori $S^{1} \times \gp$ are not isotopic in $X_{L}$ either. 

The symplectic case of problem (P) for the link surgery manifold $X_{L}$ is transformed this way into the problem of finding homologous, nonisotopic curves in $S^{3} \setminus \nu L$ transverse to the fibration. 

\textbf{Example: Braiding construction 1.} The first example we are going to consider is when $L$ is the unknot $K$. In this case, a curve transverse to the trivial fibration 
of $S^{3} \setminus \nu K$ is just a ordinary braid (closing to a knot $\gp$) having the unknot $K$ as axis. The number of strands of the braid determines the linking number $\mbox{lk}(\gp,K)$ and therefore the homology class of $S^{1} \times \gp$.  
By choosing $X$ to be the elliptic surface $E(n)$, we can produce examples of symplectic tori representing some multiple of the homology class of the fiber $F$, identified in (\ref{glue}) with $S^{1} \times \mu(K)$. The original examples of \cite{FS3} belong to this class.

\begin{figure}[h]
\begin{center}
\includegraphics[scale=.5]{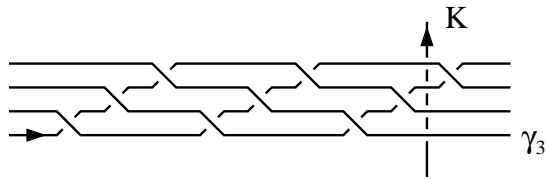}
\end{center}
\caption{The torus knot $T(3,4)$ as a braid with axis the unknot $K$.}
\label{tipiqu}
\end{figure}

A collection of braids with $q$ strands which satisfies the description above is the braid presenting the $T(p,q)$ torus knot; Figure \ref{tipiqu} illustrates this braid when $q=4,p=3$.  Here $p$ is coprime to $q$, to ensure that by closing the braid we obtain a knot.

It is clear from this construction that (at least using this approach) we cannot construct an interesting family of tori that represent the homology class of $S^{1} \times \mu(K)$ i.e. the class of the fiber in $E(n)$, for the simple fact that the only braid with one strand is the trivial one. This prevents us, in particular, from constructing nonisotopic symplectic representatives of a primitive homology class.

\begin{figure}[h]
\begin{center}
\includegraphics[scale=.5]{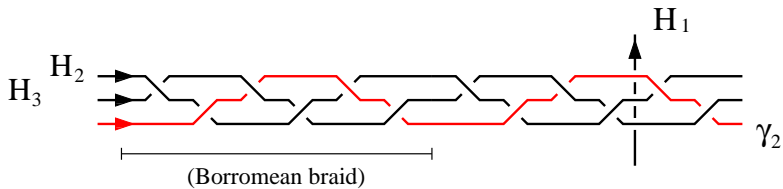}
\end{center}
\caption{Iterated Borromean rings with their axis.}
\label{itborr}
\end{figure}

\textbf{Example: Braiding construction 2.} There is a way, presented in \cite{V1}, to circumvent the difficulty of representing the fiber class, at least for $n \geq 3$, based on the fact that we can represent $E(n)$ as link surgery manifold
of the $3$-component link $L = H_{1} \cup H_{2} \cup H_{3}$ given by the Hopf link $H_{1} \cup H_{2}$ plus a third component $H_{3}$ obtained by pushing a copy of $H_{2}$ along its $0$-framing (a necklace with two rings). The observation that $S^{1} \times (S^{3} \setminus \nu L) = T^{2} \times (D^{2} \setminus \nu \{p_{1},p_{2}\})$ shows that, by choosing $X_{1} = E(n-2)$, $X_{2} = X_{3} = E(1)$, and as $\Sigma$ the spanning surface $(D^{2} \setminus \nu \{p_{1},p_{2}\})$ of $H_{1}$ pierced by $H_{2},H_{3}$ we have $X_{L} = E(n)$, for $n$ at least $3$. At this point we consider the family $\gp$ presented in Figure \ref{itborr}. 

To motivate such a choice, observe that each ``braided block" is given by the braid producing the Borromean rings. The linking matrix of $L \cup \gamma_{p}$ therefore does not depend on $p$ and this gives
$$ [\gamma_{p}] = \sum_{i=1}^{3}lk(\gamma_{p},H_{i}) [\mu(H_{i})]
= [\mu(H_{1})].  $$ The tori $S^{1} \times \gp$ represent therefore the homology class of $S^{1} \times \mu(H_{1})$, identified in (\ref{glue}) with the fiber of $E(n)$. Although the $\gp$'s are homologous,
the nontriviality of the Borromean braid translates into the fact that they are not isotopic in $S^{3} \setminus \nu L$, and we can reasonably expect that the homologous tori $S^{1} \times \gp$ are not isotopic in $E(n)$. 

Other constructions based on braiding, that produce symplectic representatives in various other classes of elliptic surfaces, are contained in \cite{EtP1} and \cite{EtP2}.

To complete the dictionary in the introduction, we have to discuss lagrangian tori. We start from the following observation. Consider a curve $\gamma \subset  S^{3} \setminus \nu L$, 
lying on a fiber $\Sigma \subset S^{3} \setminus \nu L$. After Dehn surgery, $\gamma$ is contained in ${\hat \Sigma} \subset N_{L}$ and satisfies therefore $\alpha|_{\gamma} = 0$. The torus $S^{1} \times \gamma \subset S^{1} \times N_{L}$ is therefore lagrangian, as 
$$  (dt \wedge \alpha +
\epsilon \beta)|_{S^{1} \times \gamma} = 0,  $$ and so is its image in $X_{L}$ after the gluing. For example, it is clear from this that the ``rim torus" $S^{1} \times \lambda(K) = S^{1} \times \partial D^{2} \subset E(n)$, obtained from the trivial link surgery construction with the unknot, is lagrangian. To produce infinitely many lagrangian tori we have to look for a family $\gp$ of 
curves lying on a fiber of $S^{3} \setminus \nu L$ that are not isotopic. 
It is rather clear that in order to obtain such a family we need to use a knot or a link with a fiber having a richer topology than a punctured sphere. We will discuss two examples. 

\begin{figure}[h]
\begin{center}
\includegraphics[scale=.5]{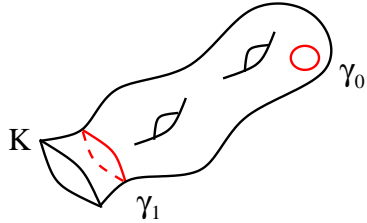}
\end{center}
\caption{A fiber with two homologous nonisotopic curves.}
\label{general}
\end{figure}

\textbf{Example: Toy model.} The first example we consider is a finite toy model (that has anyhow, as we will see later, something to teach us), that works for any nontrivial fibered knot $K$; in the fiber $\Sigma$ we consider the two homologous curves $\gamma_{0}$, a contractible curve, and a copy of the longitude of $K$ itself, that we denote as $\gamma_{1}$. These two curves are illustrated in Figure \ref{general}. 

As long as $K$ is nontrivial, $\gamma_{0}$ and $\gamma_{1}$ are not isotopic, as they are respectively unknotted and knotted in $S^{3}$. 

\begin{figure}[h]
\begin{center}
\includegraphics[scale=.5]{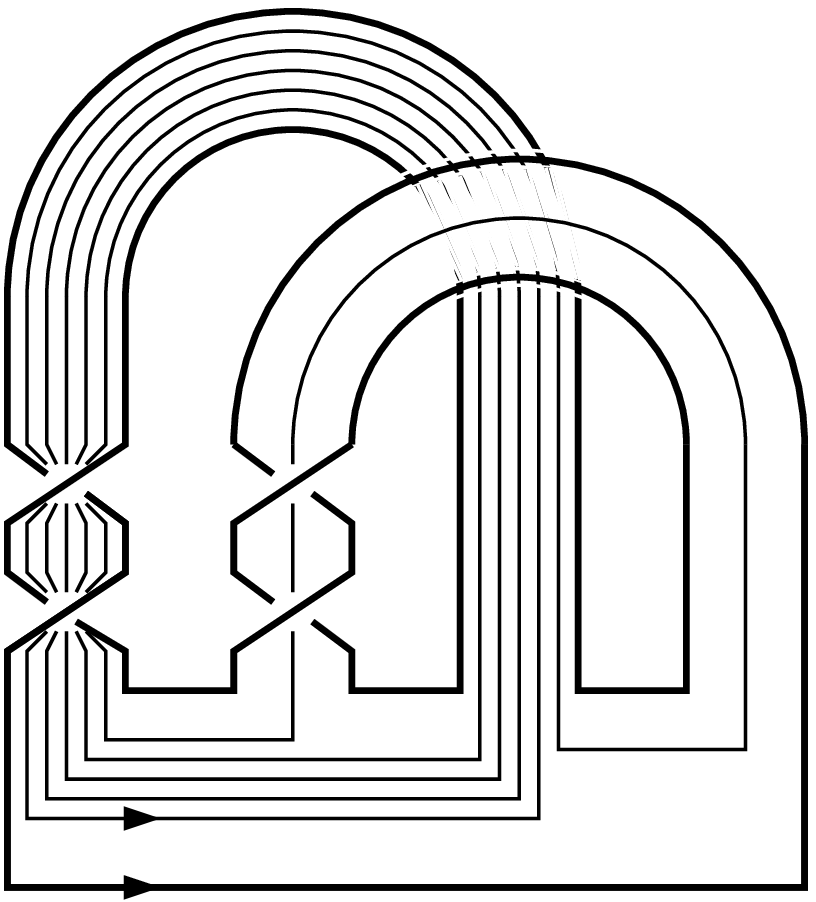}
\hspace*{2cm} \includegraphics[scale=.5]{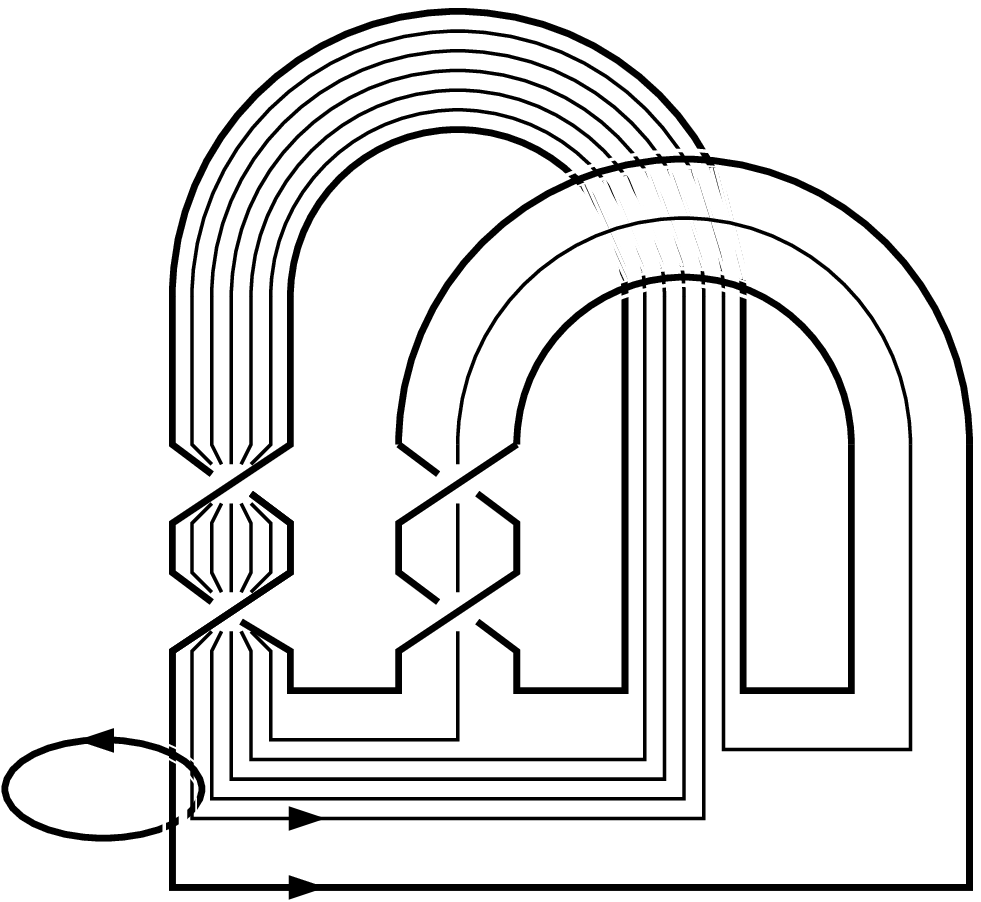}
\end{center}
\caption{The curves $\gp$ on a fiber of the trefoil $K$ (l.h.s.) and with linking number $1$ with the meridian $M$ (r.h.s.); here, $p=5$.}
\label{torusfive}
\end{figure}

\textbf{Example: Torus knots on the trefoil fiber.} The second example we consider, taken from \cite{V2}, is an infinite family $\{ \gp \}$ of curves in the fiber of the trefoil $K$ (or, by simple modification, the figure eight knot or any knot containing one of these knots as Murasugi-summands). The family is presented in the l.h.s. of Figure \ref{torusfive}; it is not too hard to verify that, as a knot in $S^{3}$, $\gp$ is the torus knot $T(p,p+1)$.

Note that in both examples above the linking number of $\gp$ and $K$ is $0$, as $\gp$ lies in the spanning surface of $K$. The lagrangian tori $S^{1} \times \gp$ of $X_{K}$ are therefore nullhomologous. To obtain lagrangian tori with nontrivial homology class we can consider a $2$-component link given by $K$ and a copy $M$ of its meridian; such a link admits an obvious fibration, the fiber being the spanning surface of $K$ pierced once by $M$. (Note that this is not the spanning surface of the link; its dual cohomology class is in fact $(1,0) \in H^{1}(S^{3} \setminus \nu L, \Z)$.) By suitably positioning the curves $\gp$ above, we can get $lk(\gp,M) = 1$ for the toy model, and $lk(\gp,M) = q$ for the trefoil case, where $q$ is any positive number and we require $p \geq q$. The r.h.s. of Figure \ref{torusfive}  illustrates such a case, with $q=1$.

In the latter construction the lagrangian tori $S^{1} \times \gp$ satisfy 
$$ [S^{1} \times \gp] = \mbox{lk}(\gp,M)[S^{1} \times \mu(M)] $$ in $X_{L}$, where $L = K \cup M$ and the gluing conditions of (\ref{glue}) determined by the fiber spanning $K$ pierced once by $M$. It is not difficult to verify that, by choosing $X_{1} = X_{2} = E(1)$, the manifold $X_{L}$ actually coincides with the knot surgery manifold $E(2)_{K}$. 

We can return now to the symplectic case to try to take advantage of the constructions of lagrangian tori. First, we want to point out that whenever we can construct a homologically essential lagrangian torus, we can make it symplectic by a perturbation of the symplectic form (see \cite{G}). Therefore, the constructions of nonisotopic lagrangian tori immediately provide examples of nonisotopic symplectic tori whenever these tori are essential (the proof of nonisotopy depends in fact only on \textit{smooth} properties of the tori). Second, we want to take advantage of the fact, observed above, that we can find nullhomologous curves, on a fiber $\Sigma$, that are nonisotopic. If we can ``sum" these curves to a fixed transverse curve, we can hope that the resulting symplectic tori inherit the nonisotopy. This is the key point of the next example.

\begin{figure}[h]
\begin{center}
\includegraphics[scale=.5]{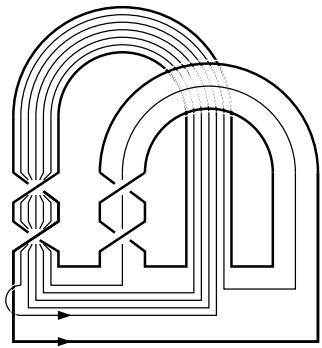}
\end{center}
\caption{Circle sum of a copy of the meridian and a nullhomologous curve; here, $p=5$.}
\label{gamma}
\end{figure}

\textbf{Example: Cabling construction.} We can apply the idea mentioned above in the following way (see \cite{V3}, \cite{EtP3}). Take again a fibered knot $K$ and consider a copy of the longitude and a contractible curve on a fiber, as in the toy model discussed above. We denote by $\gp$  the $(p,1)$-cable  of  $K$. It is easy to see that such cables, up to isotopy,  are transverse to the fibration of $S^{3} \setminus \nu K$, so that, in $X_{K}$, the tori $S^{1} \times \gp$ are symplectic.  Their homology class is as usual determined by $\mbox{lk}(K,\gp) =1$, so that $[S^{1} \times \gp] = [S^{1} \times \mu(K)] = [F]$. We can say that different values of $p$ correspond to the different number of times we have \textit{circle summed} the longitude (the $(1,0)$-cable) to the meridian (the $(0,1)$-cable). Note that, also in this case, this construction makes sense only for a nontrivial $K$. In a more general vein, we can consider curves, transverse to the fibration, obtained by circle summing a copy (or $q$ copies) of the meridian to nullhomologous, nonisotopic curves as the ones presented in Figure \ref{torusfive}. An example of that is given in Figure \ref{gamma}. This construction can be applied, for example, when $X_{1} = E(1)$; in this case the tori represent 
the homology class $[F]$ of $E(1)_{K}$. In particular, this is a primitive homology class.

\section{Detecting nonisotopy}\label{noniso}

The common feature of the constructions of the previous section is a family of homologous curves, transverse to, or  contained in, a fiber of $S^{3} \setminus \nu L$, and that are (or at least appear) nonisotopic. In fact, the proof of nonisotopy in $S^{3} \setminus \nu L$ for the examples of Section \ref{const} reduces to an obvious observation (as in the case where the $\gp$'s are different torus knots) or to the computation of the Alexander polynomial of the link $L \cup \gp$, that depends on the isotopy class of $\gp$. However, the fact that the $\gp$'s are not isotopic does not imply that the tori $S^{1} \times \gp \subset X_{L}$ are not isotopic, and this requires us to use a direct method to decide this question. 

Several methods can be used. A useful approach to this question, suggested in \cite{V1}, is based on the observation that, if two tori $T_{p_{i}}$, $i=1,2$ are isotopic, there exist a diffeomorphism of their exteriors $X_{L} \setminus \nu T_{p_{i}}$. Such a diffeomorphism implies, on its own, that there exist a diffeomorphism of the manifolds that we obtain by fiber summing $E(1)$ identifying $T_{p_{i}}$ and the fiber $F$ of $E(1)$, i.e. $X_{L} \#_{T_{p_{i}} = F} E(1)$. A word of caution is needed here. In general, the definition of fiber sum requires a choice of a diffeomorphism between $F$ and $T_{p_{i}}$, and a lifting of such a diffeomorphism to the (topologically trivial) circle bundle over these tori. In general, such choices affect the smooth type of the resulting manifold. In particular, it is \textit{a priori} impossible to establish if an isotopy between $T_{p_{1}}$ and $T_{p_{2}}$ preserves the choice of a homology basis for these tori (as we are not assuming the stronger condition that the marked tori are images of \textit{isotopic embeddings}), or their framing. However, the choice of $E(1)$ makes the situation much simpler, as any orientation preserving diffeomorphism of $\partial \nu F$ extends to an orientation preserving diffeomorphism of $E(1) \setminus \nu F$ (see e.g. \cite{GS}), so that the smooth type of $X_{L} \#_{T_{p_{i}} = F} E(1)$ is in fact unaffected by the choice of the gluing map, in particular it depends only on the isotopy type of $T_{p_{i}}$. More is actually true, namely it depends only on the diffeomorphism type of the pair $(X_{L},T_{p_{i}})$, or even of the exteriors; the smooth type of $X_{L} \#_{T_{p_{i}} = F} E(1)$ will therefore distinguish the tori $T_{p_{i}}$ up to that equivalence. Therefore we can try to detect the isotopy type of $T_{p_{i}}$ by computing the Seiberg-Witten invariants of $X_{L} \#_{T_{p_{i}} = F} E(1)$ (i.e. the relative Seiberg-Witten invariants of the pair $(X_{L},T_{p_{i}})$). This technique can be applied (with various degrees of difficulty) to all the constructions of Section \ref{const}, but 
to make the discussion more concrete, we are going to consider two explicit cases, namely the family of symplectic tori arising from the family of transverse curves $\gp$ of Figure \ref{gamma}, and the family of lagrangian tori arising from the family $\gp$ in the r.h.s. of Figure \ref{torusfive},  for the case where $q=1$. Remember than in the first case the manifold $X_{L}$ is the knot surgery manifold obtained from (say) $E(1)$ and the trefoil knot $K$, while in the second case $X_{L}$ is obtained by using two copies of $E(1)$ and the link $L = K \cup M$. 

Let's proceed with the first case: the manifold $X_{L} \#_{S^{1} \times \gp = F} E(1)$, then, is itself a link surgery manifold, as it is obtained as in (\ref{glue}) from the $2$-component link $K \cup \gp$ summing copies of $E(1)$ 
along $S^{1} \times \mu(K)$ and $S^{1} \times \gp$ (as remarked above, the choice of the gluing diffeomorphism is immaterial). We can compute the Seiberg-Witten polynomial as in (\ref{alex}), (using the fact that $SW_{(E(1),F)} = 1$) to get \begin{equation} \label{mypol} SW_{X_{L} \#_{S^{1} \times \gp = F} E(1)} = \Delta_{K \cup \gp} (t^{2}_{1},t^{2}_{2}). \end{equation} The task of computing the Alexander polynomial of the link $K \cup \gp$ of Figure \ref{gamma}, or of other classes of links, appears discouraging. Moreover, to prove that the Seiberg-Witten polynomial of (\ref{mypol}) effectively distinguishes the smooth type of $X_{L} \#_{S^{1} \times \gp = F}E(1)$ for different values of $p$, we need to show that there is no automorphism of the second homology group transforming one polynomial into the other. Those are not simple problems, at least in the general case. There is however a simple shortcut that we can take in order to obtain, essentially, the same kind of answer to problem (P), namely exhibit infinitely many nonequivalent tori. In fact, if we are able to show that the number $\beta_{p}$ of Seiberg-Witten basic classes of $X_{L} \#_{S^{1} \times \gp = F} E(1)$ satifies $\lim_{p} \beta_{p} = + \infty$, we are granted that infinitely many among the $X_{L} \#_{S^{1} \times \gp = F} E(1)$ are 
non-diffeomorphic, as $\beta_{p} < + \infty$. To prove this is a rather simpler task: in fact this amounts to showing that the number of nonzero terms of  $\Delta_{K \cup \gp}$ satisfies that property. We can simplify this task even further: using Torres formula we can observe that $$  \Delta_{K \cup \gp}(1,t_{2}) = \frac{t_{2}^{lk(K,\gp)}-1}{t_{2}-1}\Delta_{\gp}(t_{2}) = \Delta_{\gp}(t_{2}) . $$ As a knot in $S^{3}$, 
$\gp$ is the torus knot $T(p,p+1)$: using the formula for the Alexander polynomial of torus knots, we obtain that $\beta_{p}$ is bounded below by the number of nonzero terms in the Alexander polynomial of $T(p,p+1)$, namely
$$ \Delta_{T(p,p+1)}(t) = \frac{(t^{p(p+1)}-1)(t-1)}{(t^{p}-1)(t^{p+1}-1)}. $$ But it is at this point an easy exercise (see \cite{V2}) to verify that such number grows with no bound with $p$. 

The second case, for the lagrangian tori arising from Figure \ref{torusfive}, follows almost \textit{verbatim}: here, the manifold
$X_{L} \#_{S^{1} \times \gp = F} E(1)$ is itself a link surgery manifold obtained as in (\ref{glue}) from the $3$-component link $K \cup M \cup \gp$ summing copies of $E(1)$ along $S^{1} \times \mu(K)$, $S^{1} \times \lambda(M)$ and $S^{1} \times \gp$. Equation (\ref{alex}) gives now
$$ SW_{X_{L} \#_{S^{1} \times \gp = F} E(1)} = \Delta_{K \cup M \cup \gp} (t^{2}_{1},t^{2}_{2},t_{3}^{2}) $$ and, again by application of Torres formula, the number of basic classes of that manifold
is bounded below by the number of nonzero terms in the polynomial $$ \Delta_{K \cup M \cup \gp} (1,t_{2},t_{3}) = 
(t_{2}^{lk(K,M)}t_{3}^{lk(K,\gp)}-1) \Delta_{M \cup \gp}(t_{2},t_{3}) = (t_{2}-1) \Delta_{T(p,p+1)}(t_{3}). $$

To obtain the last equality we have used the fact that the pair $M \cup \gp$ is given by the torus knot $T(p,p+1)$ and its meridian $M$ , so that $\Delta_{M \cup \gp}(t_{2},t_{3}) = \Delta_{\gp}(t_{3})$. It is clear at this point that also in this case the number of basic classes is bounded below by (twice) the number of nonzero terms of the Alexander polynomial of $T(p,p+1)$. 

We want to point out the existence of another way to use Seiberg-Witten theory for detecting nonisotopy of lagrangian tori, namely the lagrangian framing defect. Such invariant, whose definition depends \textit{a priori} on symplectic properties of a lagrangian torus, can be shown under favorable circumstances to be a smooth invariant, relating it to Seiberg-Witten invariants. A detailed description of this invariant appears in \cite{FS4}. 

Finally, we discuss the simplest invariant that can detect nonisotopy of tori (or other submanifolds), namely the fundamental group of the exterior. In this case, different fundamental groups imply even the absence of pair homeomorphism.
For some of the constructions of Section \ref{const} (for example, the braiding constructions of \cite{FS3} and \cite{V1}) this invariant gives no information, while it is sufficient in other cases (for example, the cabling construction, see \cite{PV}) to distinguish the tori even up to pair homeomorphism. We will discuss its use in the following section.

\section{Higher genus examples}\label{last} 
All the examples discussed in the previous sections concerned the case of tori.
Examples of nonisotopic symplectic surfaces of genus greater than $1$ have been somewhat elusive (at least in our setup; otherwise, examples in non-simply connected manifolds have been presented in \cite{S}, and singular examples in \cite{ADK}).
One of the reason for this is the difficulty in using Seiberg-Witten poynomials, as in Section \ref{noniso}, to detect potential higher genus examples built by ``doubling" the braiding constructions considered in Section \ref{const}. In this section we are going to discuss how we can instead double the cabling construction and retain enough information from the fundamental group of the exterior to be able to distinguish higher genus examples. For sake of simplicity we will consider just the case of genus $2$, but the general case follows along the same lines. For details, see \cite{PV}. Start by considering the knot surgery manifold $E(2)_{K}$ obtained from a nontrivial fibered knot $K$ choosing $X_{1} = E(2)$. In such a manifold we can identify a 
symplectic surface of genus $g(K)$ and self-intersection $-2$ obtained by capping off the spanning surface $\Sigma$ of the knot with a $-2$-disk section in $E(2) \setminus \nu F$, according to the gluing prescription of (\ref{glue}). Summing such a surface with a copy of the elliptic fiber, and suitably resolving, we get a symplectic surface $\Sigma_{g}$ of genus $g = g(K)+1$ and self-intersection $0$. Following \cite{P} we can define a symplectic, simply connected $4$-manifold by fiber summing two copies of $E(2)_{K}$ along $\Sigma_{g}$, i.e. $$ D_{K} = (E(2)_{K} \setminus \nu \Sigma_{g}) \cup (E(2)_{K} \setminus \nu \Sigma_{g}) = E(2)_{K} \#_{\Sigma_{g}} E(2)_{K},  $$ identifying the two boundary components $\partial \nu \Sigma_{g} = S^{1} \times \Sigma_{g}$  with a diffeomorphism that is the identity on $\Sigma_{g}$ and reverses the orientation on $S^{1}$. In $E(2)_{K}$ we have the family of homologous symplectic tori $S^{1} \times \gp$ that we introduced in Section \ref{const} through the cabling construction. These tori intersect the surface $\Sigma_{g}$ in a single transverse point, inherited from the intersection point between $\gp$ and the spanning surface of $K$. It follows that in $D_{K}$ the two copies of $S^{1} \times \gp$ located in each half give, by connected sum, a genus $2$ surface $\Xi_{p}$. 

We can attempt to distinguish these surfaces using the fundamental group of the exteriors $D_{K} \setminus \nu \Xi_{p}$. A rather lengthy exercise in Van Kampen theorem (see \cite{PV}) shows that this group is given by \begin{equation} \label{group} \pi_{1}(D_{K} \setminus \nu \Xi_{p}) = \frac{\pi_{1}(S^{3} \setminus \nu K)}{\mu(K)\lambda(K)^{p}} = \pi_{1}(S^{3}_{1/p}(K)) \end{equation} i.e. the fundamental group of the homology sphere obtained by applying Dehn surgery to $S^{3}$ along $K$ with slope $1/p$. The problem
of distinguishing the surfaces $\Xi_{p}$ boils down therefore to showing that, for a given knot $K$, infinitely many among these homology spheres are distinguished by means of their fundamental group. Before commenting on that, note that if $K$ is the unknot, this fundamental group is always trivial: this is exactly what we expect, as it follows immediately from the construction of Section \ref{const} that in that case  the tori $S^{1} \times \gp$ are isotopic (the $(p,1)$-cable of the unknot is isotopic to the meridian!). 

Let us return to the problem of distinguishing up to isomorphism the groups of (\ref{group}). We will consider here the case where $K$ is a fibered hyperbolic knot. In this case, with a finite number of exception, the homology spheres $S^{3}_{1/p}(K)$ are hyperbolic $3$-manifolds. Mostow rigidity implies in that case that these manifolds are homeomorphic if and only if their fundamental groups coincide. Therefore, we are left with the problem of showing that infinitely many of the $S^{3}_{1/p}(K)$ are not homeomorphic. But at this point we can appeal to Theorem 1 of \cite{BHW}, that states that hyperbolic fillings of a hyperbolic knot exterior can be homeomorphic only if the knot is amphicheiral and the slopes are opposite. Restricting ourselves to positive $p$ for sake of simplicity, we deduce that for all but finitely many choices of $p$ the fundamental group of (\ref{group}) distinguishes (up to homeomorphism the exteriors, hence isotopy) the genus $2$ surfaces $\Xi_{p}$. An analogous result can be proven for all nontrivial fibered knots.

\end{document}